\documentclass[11pt]{conm-p-l}
\usepackage{cite}

\newtheorem{theorem}{Theorem}[section]
\newtheorem{lemma}[theorem]{Lemma}

\theoremstyle{definition}
\newtheorem{definition}[theorem]{Definition}

\theoremstyle{remark}
\newtheorem{remark}[theorem]{Remark}

\numberwithin{equation}{section}

\newtheorem{corollary}[theorem]{Corollary}

\def\Q{\mathbb{Q}}
\def\R{\mathbb{R}}
\def\C{\mathbb{C}}

\def\F{\mathbb{F}}
\def\1{\textbf{1}}

\begin{document}

\title[Non-commutative Binomials]{Some binomial formulas for non-commuting operators}

\author{Peter Kuchment}
\address{Mathematics Department, Texas A\&M University, College Station, TX 77845}
\email{kuchment@math.tamu.edu}
\thanks{P. K. expresses his gratitude to the support from DMS NSF Grant \#1211463.}
\author{Sergey Lvin}
\address{Mathematics and Statistic Department, University of Maine, Orono, ME 04469}
\email{lvin@math.umaine.edu}
\dedicatory{Dedicated to the memory of our beloved teacher, colleague, and co-author Selim Krein}

\date{}

\begin{abstract}Let $D$ and $U$ be linear operators in a vector space (or more generally, elements of an associative algebra with a unit). We establish
binomial-type identities for $D$ and $U$ assuming that either their commutator $[D,U]$ or the second commutator $[D,[D,U]]$ is
proportional to $U$.

Operators $D=d/dx$ (differentiation) and $U$- multiplication by $e^{\lambda
x}$ or by $\sin \lambda x$ are basic examples, for which some of these relations appeared unexpectedly as byproducts of an authors' medical imaging
research \cite{KL89,KL90,KL13,Kuchbook}.
\end{abstract}

\maketitle
\section*{Introduction}

While working on range conditions for a Radon type transform arising in emission medical imaging, the authors \cite{KL89,KL90} (see also \cite{KL13}) discovered that one of their theorems was equivalent to an infinite series of puzzling nonlinear combinatorial-differential identities for the classical exponential, linear, and some trigonometric functions (sic!). Here are the examples:

\begin{itemize}
\item For any non-negative integer $n$ and $u=e^{\lambda x}$, one has
\begin{equation}\label{E:oldexp}
\sum\limits_{k=0}^n \dbinom{n}{k}
\left[\left(\frac{d}{dx}-u\right)\circ \left(\frac{d}{dx}-u+\lambda\right)\circ\dots\circ\left(\frac{d}{dx}-u+(k-1)\lambda\right)\right]u^{n-k}=0.
\end{equation}
\item For $u=\sin\lambda x$, similar identities hold for any odd natural $n$:
\begin{equation}\label{E:oldsin}
\sum\limits_{k=0}^n \dbinom{n}{k}
\left[\left(\frac{d}{dx}-u\right)\circ \left(\frac{d}{dx}-u+i\lambda\right)\circ\dots\circ\left(\frac{d}{dx}-u+i(k-1)\lambda\right)\right] u^{n-k}=0.
\end{equation}
\end{itemize}
In order to avoid quite possible misunderstanding, let us explain briefly the meaning of various terms of these identities.
\begin{remark}
The factors in the $\circ$-products are understood as operators on smooth functions on the line. In particular, $u$ there means the multiplication operator by the function $u$ (exponential or sine). The $\circ$ means the composition of operators, and the order of factors is important, due to them non-commuting. On the other hand, $u^{n-k}$ at the end is considered as a function, to which the operator $[...]$ is applied. If one tries to understand $u^{n-k}$ also as an operator and thus considers $\sum \dbinom{n}{k}[\cdots]\circ u^{n-k}$, the resulting \textbf{operator} is NOT identically equal to zero. It needs to be applied to the function identically equal to $1$ to preserve the identity.
This understanding will be important throughout the text.

It is also why, to avoid misinterpretation, \textbf{we use later the notation $\widetilde{f}$ for the operator of multiplication by the function $f$}.
\end{remark}

In the paper, we significantly extend these results, as well as generalize them to a wider algebraic situation. Namely, the setting in which we obtained these results before was of a commutative algebra (where $u$ belongs to) with differentiation $D$ and $u$ satisfying ``differential equations'' $Du=\lambda u$ or $D^2u=\lambda^2 u$. Now we show that the results have generalization to elements $D$ and $u$ of any associative algebra with a unit, with appropriate conditions on their first and second commutators.

 When these identities appeared in medical imaging, they have attracted quite a lot of attention, especially after discovering their relations to unusual Hartogs type analytic continuation theorems in several complex variables (\cite{AEK,Oktem1,Oktem2,Tum}). It is fair to notice that, in spite of a variety of different proofs and generalizations of the identities available, the authors still feel that they do not have good understanding of the origin of such formulas.


here is the structure of the paper: The main notions are introduced and results stated in Section  \ref{S:results}. Section \ref{S:funct} is devoted to the particular cases of elementary functions. The proofs are delegated to Section \ref{S:proofs}, followed by a final remarks section.

\section{Formulation of main results}\label{S:results}

Let $D$ and $U$ be elements of an associative algebra $A$ over a field $\Q$ with identity $I$ (for instance, the algebra of all linear operators in a vector
space $\mathbb{F}$).

Let us introduce,  lead by (\ref{E:oldexp})-(\ref{E:oldsin}), an $n$th order binomial-type
combination of $D$ and $U$
\begin{equation}B(n,\lambda,U,D):=\sum\limits\limits_{k=0}^{n}\dbinom{n}{k}\left(
\prod\limits_{j=0}^{k-1}\left( D-U+j\lambda I\right) \right) U^{n-k},
\label{Bdef}
\end{equation}
where $n\geq 0$ is an
integer, $\lambda \in\Q$, and $\dbinom{n}{k}$ is the binomial coefficient. When $k=0,$ the product is understood as $I$.

\begin{center}
\textbf{Due to non-commutativity of $A$, we will adhere to the following agreement:
The products $\prod_j \cdots$ are understood in the order of the index $j$ increasing from left to right.}
\end{center}

For instance, $B(0,\lambda ,U,D)=I$, $B(1,\lambda ,U,D)=U+(D-U)=D$, $B(2,\lambda ,U,D)=U^{2}+2(D-U)U+(D-U)(D-U+\lambda I)$, etc., with formulas getting more complex with $n$ increasing.

\subsection{First order commutator}\label{SS:single}

In the following theorems we make certain assumptions about the commutator $[D,U]=DU-UD$.

The following rather surprising result holds:
\begin{theorem}\label{T:No U}
Suppose that $\left[ D,U\right] =\lambda U$. Then $B(n,\lambda ,U,D)$
does not depend on $U$. Moreover,
\begin{equation}
B(n,\lambda ,U,D)=B(n,\lambda ,0,D)=\prod\limits_{j=0}^{n-1}\left(
D+j\lambda I\right) .
\end{equation}
\end{theorem}

\begin{remark}\indent
\begin{itemize}
\item If $D$ and $U$ commute ($\lambda =0$), then the theorem states that $$B(n,0,U,D)=D^{n},$$ which is an obvious consequence of the standard binomial
formula. Indeed, when $D$ and $U$ commute,
$$B(n,0,U,D)=\sum\limits_{k=0}^{n}\dbinom{n}{k}\left( D-U\right)
^{k}U^{n-k}=(D-U+U)^{n}=D^{n}.$$
\item The equality $[D,U]=\lambda U$ is homogeneous of degree one with respect to $(D,\lambda)$. Homogeneity is not obvious for the originally defined $B(n,\lambda,U,D)$, however the statement of the theorem implies that the homogeneity does hold:
\begin{equation}
B(n,\lambda ,U,D)=\lambda^n B(n,1 ,U,D/\lambda),
\end{equation} This shows that essentially the study boils down to only the cases when $\lambda=0$ (considered above) and $\lambda=1$,
which simplifies considerations.
\end{itemize}
\end{remark}
Here are some immediate consequences of the Theorem:

\begin{corollary}\label{C:U+W} Suppose $\left[D,U\right] =\lambda U$. Then,
\begin{enumerate}
\item
If $V\in A$ and for some $j\in \{0, ..., n-1\}$, one has $(D+j\lambda I)V=0$, then
\begin{equation}
\left(B(n,\lambda
,U,D)\right)V=0.
\end{equation}

If $A$ is the algebra of linear operators on a vector space, the equivalent reformulation, under the same assumptions, is:
$$
\emph{Ker} (D+j\lambda I)\subset \emph{Ker} (B(n,\lambda
,U,D)).
$$
\item If $[V,W]=0$ and $[D,W]=\lambda W$, then $B(n,\lambda ,V+W,D)$ does not depend on $W$: $$B(n,\lambda ,V+W,D)=B(n,\lambda ,V,D).$$
Indeed, just substitute $D$ in the Theorem with $D-V$.
\end{enumerate}
\end{corollary}

For what follows, it is interesting to understand what happens with  $B(n,\lambda ,U,D)$ when $\left[ D,U\right] =-\lambda U$ (notice the wrong sign in the commutation relation, and thus Theorem \ref{T:No U} describes $B(n,-\lambda ,U,D)$, rather than $B(n,\lambda ,U,D)$). Things get more complicated here.

In the following results we denote by $\mathbb{F}_{0}$ the set of all $V\in A$ is such that $DV=0$. We use the standard notation $(2m-1)!!$ for $1\cdot 3\cdot
...\cdot (2m-1).$

\begin{theorem}\label{T:wrongsign} Let $\left[ D,U\right] =-\lambda U$, then
\begin{itemize}
\item $B(n,\lambda ,U,D)|_{\mathbb{F}_{0}}=0$ for any odd $n$,
\item $B(n,\lambda ,U,D)|_{\mathbb{F}_{0}}=(n-1)!!(-2\lambda U)^{n/2}|_{\mathbb{F}_{0}}$ for even $n>0$, and
\item $(2D+\lambda nI)B(n,\lambda ,U,D)|_{\mathbb{F}_{0}}=0$ for all $n\geq
0$.
\end{itemize}
\end{theorem}
The reason why we interested in this case will be clear in the proof of Theorem \ref{T:2ndU} below.
\subsection{Second order commutators}\label{SS:double}

Here we will be interested in using conditions on the second order commutators $[D,[D,U]]$ and $[U,[D,U]]$. Here, the condition $[D,U]=\lambda U$ is replaced with $[D,[D,U]]=\lambda ^{2}U$ (notice that we preserve the $(D,\lambda)$-homogeneity)\footnote{For elementary functions example (see the next section), this means switching from exponential to trigonometric functions.}. We are interested in behavior of the polynomials $B(n,\lambda ,U,D)$ in this situation.

\begin{theorem}\label{T:2ndU} Suppose $[D,[D,U]]=\lambda ^{2}U$ and $[U,[D,U]]=0$. Then
\begin{itemize}
\item $B(n,\lambda ,U,D)|_{\mathbb{F}_{0}}=0$ for all odd $n$,
\item $B(n,\lambda ,U,D)|_{\mathbb{F}_{0}}=(n-1)!!([D,U]-\lambda U)^{n/2}|_{\mathbb{F}_{0}}$  for
all even $n>0$, and
\item $(2D+\lambda nI)B(n,\lambda ,U,D)|_{\mathbb{F}_{0}}=0$ for all $n\geq
0$.
\end{itemize}
\end{theorem}

\begin{remark} Theorem \ref{T:2ndU} (unlike Theorem \ref{T:No U}) remains nontrivial even when $\lambda =0.$
\end{remark}

\section{Differential identities for some elementary functions}\label{S:funct}
Here we apply the above results to the algebra of linear operators acting in the vector space $\F:=C^\infty(\R)$
of all smooth functions on the real line.
 \begin{definition}\indent
 \begin{itemize}
 \item We \textbf{denote by $\widetilde{f}$ the operator of multiplication by such a function $f(x)$}.
 \item We use the notation $\1$ for the function that is identically equal to $1$.
 \item We also denote $D:=d/dx$.
 \end{itemize}
 \end{definition}

\subsection{Exponential functions}\label{SS:exp}
It is clear that $e^{-j\lambda x}$ belongs the kernel of $D+j\lambda I$ and, in particular, $1$ is an element in the kernel of $D$.

Now one clearly has
$$[D,\widetilde{e^{\pm \lambda x}}]=\pm\lambda \widetilde{e^{\pm \lambda x}}.$$
Thus, the results of Section \ref{SS:single} apply to produce the following formulas:
\begin{theorem}\label{T:newexp}
 \begin{equation}
B(n,\lambda ,d/dx,\widetilde{e^{\lambda x}})e^{-j\lambda x}=0\mbox{ for all }n>0\mbox{
and }0\leq j\leq n-1,
\end{equation}
\begin{equation}
B(n,\lambda ,d/dx,\widetilde{e^{-\lambda x}})\textbf{1}=0\mbox{ for all odd }n,
\end{equation}
\begin{equation}
B(n,\lambda ,d/dx,\widetilde{e^{-\lambda x}})\textbf{1}=(n-1)!!(-2\lambda )^{n/2}e^{-n\lambda
x/2}\mbox{ for all even }n>0,
\end{equation}
\begin{equation}
(2d/dx+\lambda nI)B(n,\lambda ,d/dx,\widetilde{e^{-\lambda x}})\textbf{1}=0\mbox{ for all }n>0.
\end{equation}
\end{theorem}

\subsection{Trigonometric functions}\label{SS:trig}
Since multiplications by functions commute, we have in this case the condition $[\widetilde{f},[D,\widetilde{f}]]=0$, needed in
Section \ref{SS:double}, automatically satisfied for any smooth function $f$.

Let us now check the condition $[D,[D,\widetilde{f}]]=\lambda^2 \widetilde{f}$ for the natural candidates: trigonometric and hyperbolic sine and cosine.

It is an easy computation that
when $U=\widetilde{\sin \lambda x}$ is the operator of multiplication by $\sin \lambda
x$,  then $[D,U]=\lambda \widetilde{\cos \lambda x}$, $[D,[D,\widetilde{\sin \lambda x}]]=(i\lambda )^{2}\widetilde{\sin \lambda x}$, and as we have mentioned above, $[U,[D,U]]=0$ is automatic. 

Thus, the results of section \ref{SS:double} provide the following set of differential
identities for sine functions:

\begin{theorem}\label{T:newsin}
\indent\begin{itemize}\item
\begin{equation}
B(n,i\lambda ,d/dx,\widetilde{\sin \lambda x})\textbf{1}=0\mbox{ for all odd }n,
\end{equation}
\item
\begin{equation}
B(n,i\lambda ,d/x,\widetilde{\sin \lambda x})\textbf{1}=(n-1)!!\lambda ^{n/2}e^{-in\lambda x/2}%
\mbox{ for all even }n>0,
\end{equation}
\item
\begin{equation}
(2d/dx+i\lambda nI)B(n,i\lambda ,d/x,\widetilde{\sin \lambda x})\textbf{1}=0\mbox{ for all }n>0.
\label{(c)}
\end{equation}
\end{itemize}
\end{theorem}

\begin{remark}
Similar identities hold for any solutions of the differential equation $\dfrac{%
d^{2}u}{du^{2}}=\lambda ^{2}u,$ including cosines, hyperbolic
sine and cosine, and linear functions.

In particular, here are the binomial-type identities for
linear functions:
\begin{itemize}
\item\begin{equation}
B(n,0,d/dx,\widetilde{(ax+b)})\textbf{1}=0\mbox{ for all odd }n,
\end{equation}
\item
\begin{equation}
B(n,0,d/dx,\widetilde{(ax+b)})\textbf{1}=(n-1)!!a^{n/2}\mbox{ for all even }n>0,
\end{equation}
\item
\begin{equation}
(d/dx)B(n,0,d/dx,\widetilde{(ax+b)})\textbf{1}=0\mbox{ for all }n>0.
\end{equation}%
\end{itemize}
\end{remark}

\subsection{Change of variables}

One can play with changes of variables in the formulas of Theorems \ref{T:newexp}
and \ref{T:newsin}, to get a variety of new identities. For instance, the change $x\to x^{2}/2$ and correspondingly $d/dx \to x^{-1}d/dx$ gives
\begin{equation}
B(n,\lambda ,\frac{1}{x}\frac{d}{dx},\widetilde{e^{\lambda x^{2}/2}})e^{-j\lambda x^{2}/2}=0\mbox{
for all }n>0\mbox{ and }0\leq j\leq n-1,
\end{equation}
while $x \to \ln x, d/dx \to x\,d/dx$ produces
\begin{equation}
B(n,\lambda ,x\frac{d}{dx},\widetilde{\lambda x})x^{-\lambda j}=0\mbox{ for all }n>0%
\mbox{ and }0\leq j\leq n-1.
\end{equation}

\subsection{Vector functions}\label{SS:vector}
The results easily translate to vector-valued functions. Let $\mathbb{F}$ is the space of all smooth $\C^m$-valued functions. Then the kernel of $D=d/dx$ consists of constant column
vectors $\overrightarrow{c}$. Let $A$ be an $m\times m$ constant matrix. Then
the following set of statements hold:\\
\begin{enumerate}
\item $$
B(n,\lambda ,d/dx,\widetilde{e^{\lambda x}A})e^{-j\lambda x}I=0\mbox{ for all }n>0\mbox{
and }0\leq j\leq n-1,
$$
\item
$$
B(n,\lambda ,d/dx,\widetilde{e^{-\lambda x}A})\overrightarrow{c}=\overrightarrow{0}\mbox{
for all odd }n,
$$
\item
$$
B(n,\lambda ,d/dx,\widetilde{e^{-\lambda x}A})\overrightarrow{c}=(n-1)!!(-2\lambda
)^{n/2}e^{-n\lambda x/2}A^{n/2}\overrightarrow{c}\mbox{ for all even }n>0,
$$
\item
$$
(2d/dx-\lambda nI)B(n,\lambda ,d/dx,\widetilde{e^{-\lambda x}A})\overrightarrow{c}=%
\overrightarrow{0}\mbox{ for all }n>0,
$$
\item
$$
B(n,i\lambda ,d/dx,\widetilde{\sin \lambda xA})\overrightarrow{c}=\overrightarrow{0}%
\mbox{ for all odd }n,
$$
\item
$$
B(n,i\lambda ,d/x,\widetilde{\sin \lambda xA})\overrightarrow{c}=(n-1)!!\lambda
^{n/2}e^{-in\lambda x/2}A^{n/2}\overrightarrow{c}\mbox{ for all even }n>0,
$$
\item
$$
(2d/dx+i\lambda nI)B(n,i\lambda ,d/x,\widetilde{\sin \lambda xA})\overrightarrow{c}=%
\overrightarrow{0}\mbox{ for all }n>0.
$$
\item
If $A_{1}$ and $A_{2}$ are commuting $m\times m$ matrices, then

$$
B(n,0,d/dx,\widetilde{A_{1}x+A_{2}})\overrightarrow{c}=\overrightarrow{0}\mbox{ for all
odd }n,
$$

$$
B(n,0,d/dx,\widetilde{A_{1}x+A_{2}})\overrightarrow{c}=(n-1)!!A_{1}^{n/2}\overrightarrow{%
c}\mbox{ for all even }n>0,
$$

$$
(d/dx)B(n,0,d/dx,\widetilde{A_{1}x+A_{2}})\overrightarrow{c}=\overrightarrow{0}\mbox{
for all }n>0.
$$%
\end{enumerate}
\section{Proofs}\label{S:proofs}
\subsection{Proof of Theorem \ref{T:No U}}
We start with the following lemma, which can be easily
proved by induction.

\begin{lemma}\label{L1} \indent
\begin{enumerate}
\item If $\left[ D,U\right] =\lambda U$, then $\left[ D,U^{m}\right]
=\lambda mU^{m}$ for all integer $m\geq 0$.
\item If $[D,U]=V$, $[U,V]=0$, and $[D,V]=0$, then $%
[D,U^{m}]=mVU^{m-1}$ for all natural $m.$
\end{enumerate}
\end{lemma}

We now provide a different representation for $B(n,\lambda,U,D)$.

\begin{lemma}\label{L2} For all $n>0$%
\begin{equation}
B(n,\lambda ,U,D)=\sum\limits\limits_{k=0}^{n-1}\dbinom{n-1}{k}\left(
\prod\limits_{j=0}^{k-1}(D-U+\lambda j)\right) (D+\lambda kI)U^{n-1-k}.
\label{(2)}
\end{equation}
\end{lemma}
\textbf{Proof} Using $\dbinom{n}{k}=\dbinom{n-1}{k}+\dbinom{n-1}{k-1},$ we
can rewrite $B(n,\lambda ,U,D)$ as follows:
$$
\begin{array}{c}\label{E:L2}
\sum\limits\limits_{k=0}^{n-1}\dbinom{n-1}{k}\left(
\prod\limits_{j=0}^{k-1}\left( D-U+j\lambda I\right) \right) U^{n-k}\\
+\sum\limits\limits_{k=1}^{n}\dbinom{n-1}{k-1}\left(
\prod\limits_{j=0}^{k-1}\left( D-U+j\lambda I\right) \right) U^{n-k}.
\end{array}
$$

(\textbf{We remind the reader that in the products the order is in increasing $j$ (or $k$) from the left to the right! })
Changing $k$ to $k+1$ in the second sum and combining both
sums proves the lemma. $\Box$.

\textbf{Proof of Theorem \ref{T:No U}}. By definition, $B(0,\lambda ,U,D)=I$ and $B(1,\lambda
,U,D)=D $.

Let now $n>1$ and $\left[ D,U\right] =\lambda U$.  Using Lemma \ref{L1}, we get
$$(D+\lambda kI)U^{n-1-k}=U^{n-1-k}(D+\lambda (n-1)I).$$
Substituting this into
(\ref{E:L2}), we obtain the recurrent formula
\begin{equation}
B(n,\lambda ,u,D)=B(n-1,\lambda ,u,D)(D+\lambda (n-1)I).  \label{(3)}
\end{equation}%
By induction, (\ref{(3)}) implies Theorem \ref{T:No U} for all $n$. $\Box$

\subsection{Proof of Theorem \ref{T:wrongsign}}

The following lemma (appeared in a somewhat more restricted and implicit form in \cite{KL90}) provides an important insight on the nature of $%
B(n,\lambda ,U,D)$, and the technique of its proof will be used
below. It uses an additional assumption that $U$ is invertible. It is easy to see that if $\left[ D,U\right] =\lambda U$,
then $\left[ D,U^{-1}\right] =-\lambda U^{-1},$ and the statement of Lemma \ref{L1} holds now for all
integers $m$, including $m<0$.

The following lemma provides a useful recurrent relation for $B(n,\lambda ,U,D)$ in the case when $[D,U]=-\lambda U$.

\begin{lemma}\label{L:wrongsign} If $[D,U]=-\lambda U,$ then
\begin{equation}\begin{array}{l}
B(n,\lambda ,U,D)=B(n-1,\lambda ,U,D)(D+\lambda (n-1)I)
-2(n-1)\lambda UB(n-2,\lambda ,U,D)\\+2(n-1)(n-2)\lambda ^{2}UB(n-3,\lambda
,U,D)  \label{(6)}\end{array}
\end{equation}%
for all $n>1,$ with the last term absent when $n=2.$
\end{lemma}

\textbf{Proof of Lemma \ref{L:wrongsign}} Using Lemma \ref{L1}, we get $$\begin{array}{l}(D+\lambda kI)U^{n-1-k}=U^{n-1-k}(D+\lambda
(2k-n+1)I)\\=U^{n-1-k}(D+\lambda (n-1)I-2\lambda (n-1-k)I).\end{array}$$
Substituting it
into (\ref{(6)}), we obtain
$$\begin{array}{l}
B(n,\lambda ,U,D)=\sum\limits\limits_{k=0}^{n-1}\dbinom{n-1}{k}\left(
\prod\limits_{j=0}^{k-1}(D-U+\lambda jI)\right) U^{n-1-k}(D+\lambda
(n-1)I)\\-2\lambda \sum\limits\limits_{k=0}^{n-1}\dbinom{n-1}{k}(n-1-k)\left(
\prod\limits_{j=0}^{k-1}(D-U+\lambda jI)\right) U^{n-1-k}.\end{array}$$

Here the first sum is equal to $$B(n-1,\lambda ,U,D)(D+\lambda (n-1)I),$$ in the
second sum the term with $k=n-1$ is zero, and $$\dbinom{n-1}{k}(n-1-k)=(n-1)%
\dbinom{n-2}{k}.$$ Thus,

$$\begin{array}{l}Z:=B(n,\lambda ,U,D)-B(n-1,\lambda ,U,D)(D+\lambda (n-1)I)\\=-2(n-1)\lambda \sum\limits\limits_{k=0}^{n-2}\dbinom{n-2}{k}\left(
\prod\limits_{j=0}^{k-1}(D-U+\lambda jI)\right) U^{n-1-k}.\end{array}$$
Since $(D-U+\lambda jI)U=U(D-U+\lambda (j-1)I)$, we have
$$Z=-2(n-1)\lambda U\sum\limits\limits_{k=0}^{n-2}\dbinom{n-2}{k}\left(
\prod\limits_{j=-1}^{k-2}(D-U+\lambda jI)\right) U^{n-2-k}.$$

When $j=-1$, $D-U+\lambda jI=(D-U+(k-1)\lambda I)-k\lambda I$. Thus,

$$\begin{array}{l}Z=-2(n-1)\lambda U\sum\limits\limits_{k=0}^{n-2}\dbinom{n-2}{k}\left(
\prod\limits_{j=0}^{k-1}(D-U+\lambda jI)\right) U^{n-2-k}\\+2(n-1)\lambda ^{2}U\sum\limits\limits_{k=1}^{n-2}\dbinom{n-2}{k}k\left(
\prod\limits_{j=0}^{k-2}(D-U+\lambda jI)\right) U^{n-2-k}.\end{array}$$
The first sum here is $B(n-2,\lambda ,U,D)$. Using $%
\dbinom{n-2}{k}k=(n-2)\dbinom{n-3}{k-1}$ and changing $k-1$ to $k$, we conclude that the second sum equals $(n-2)B(n-3,\lambda ,U,D)$. Therefore,
$$Z=-2(n-1)\lambda UB(n-2,\lambda ,U,D)+2(n-1)(n-2)\lambda ^{2}UB(n-3,\lambda
,U,D),$$ which coincides with (\ref{(6)}) and completes the proof. $\Box$

\textbf{Proof of Theorem \ref{T:wrongsign}}
We know that $B(0,\lambda ,U,D)=I$ and $B(1,\lambda
,U,D)=D, $ so $B(1,\lambda ,U,D)=0$ on $\mathbb{F}_{0}$ and the first statement of the Theorem holds true for $n=0$ and $n=1$. Let $[D,U]=-\lambda U$, then on $%
\mathbb{F}_{0}$, the recurrence (\ref{(6)}) is reduced to
$$\begin{array}{l}
B(n,\lambda ,U,D)=\lambda (n-1)I)B(n-1,\lambda ,U,D)
-2(n-1)\lambda UB(n-2,\lambda ,U,D)\\+2(n-1)(n-2)\lambda ^{2}UB(n-3,\lambda
,U,D).\end{array}
$$
From here, the statement of the Theorem follows by induction for all $n>0$, both odd and even.$\Box$

\subsection{Proof of Theorem \ref{T:2ndU}}
The proof will be different for $\lambda \neq 0$ and for $\lambda =0$. This is not that surprising, since the claim, if considered on elementary functions, leads to differential binomial for quite different functions: sine
functions when $\lambda \neq 0$ and linear functions when $\lambda =0$.

\textbf{The case $\lambda \neq 0$.} Let us
introduce the elements $V:=(\lambda U-[D,U])/2\lambda $ and $W:=(\lambda
U+[D,U])/2\lambda$. A direct calculation shows that $U=V+W$, and under the Theorem's assumptions $[D,V]=-\lambda V$, $[D,W]=\lambda W$, and $[V,W]=0$. From
Corollary \ref{C:U+W} we obtain that $B(n,\lambda ,U,D)=B(n,\lambda ,V+W,D)=B(n,\lambda ,V,D)$. Now one applies Theorem \ref{T:wrongsign}.

\textbf{The case $\lambda = 0$.}
The following Lemma can be easily proved by induction.

\begin{lemma}\label{L6} Suppose $[D,U]=V$, $[U,V]=0$, and $[D,V]=0$. Then $%
[D,U^{m}]=mVU^{m-1}$ for all natural $m$.
\end{lemma}
 By definition, we have  $%
B(1,0,U,D)=D$ and $B(2,0,U,D)=D^{2}+[D,U]$. Let $V:=[D,U]$. Then $%
B(1,0,U,D)=0$ and $B(2,0,U,D)=V$ on $\mathbb{F}_{0}$ so the claim is true for $%
n=1$ and $n=2$. Now, under the assumption that $[U,V]=0$ and $[D,V]=0$, we will
show that for $n>2$
\begin{equation}
B(n,0,U,D)=(n-1)VB(n-2,0,U,D)\mbox{ on }\mathbb{F}_{0}.  \label{(7)}
\end{equation}%
Then induction proves the statement for all $n$.

Let $n>2$.  According to Lemma \ref{L2},
$$
B(n,0,U,D)=\sum\limits\limits_{k=0}^{n-1}\dbinom{n-1}{k}(D-U)^{k}DU^{n-1-k},
$$
where the term with $k=n-1$ is equal to zero on $\mathbb{F}_{0}.$ Also,
according to Lemma \ref{L6}, $DU^{m}=mVU^{m-1}$ on $\mathbb{F}_{0}.$ Then on $\mathbb{F}
_{0}$ we get
$$
B(n,0,U,D)=\sum\limits\limits_{k=0}^{n-2}\dbinom{n-1}{k}%
(D-U)^{k}(n-1-k)VU^{n-2-k}.
$$
From here we obtain (\ref{(7)}) by using $\dbinom{n-1}{k}(n-1-k)=(n-1)\dbinom{n-2}{k}
$ and commuting $V$ to the left. This completes the proof. $\Box$

\section{Final remarks and conclusions}
\begin{itemize}
\item It would be interesting to figure out what can be done  under the condition of vanishing of the third commutator $[D,[D,[D,U]]]=\lambda^3 U$ (maybe plus some other restrictions). It has been checked by direct computation that the natural analog of (\ref{E:oldexp})-(\ref{E:oldsin}) for the solutions of the third order equation $D^3u=\lambda^3u$ does not hold \cite{Elaina,KL13}.
\item As we have already mentioned, we do not truly understand the origin of such identities. It looks like this issue is in the realm of the techniques of \cite{Touf}, in which we are not experts, to say the least.
\item There still might be interesting relations to SCV, as the ones to Hartogs' type theorems in \cite{AEK,Tum,Oktem1,Oktem2}. One also wonders about such higher dimension analogs of Hartogs' theorems.
\item We cannot help it providing a cute lemma used in \cite{KL90}. A version of this text used it, but we have managed to avoid this. A reader, however, could find it interesting:
\begin{lemma}\cite{KL90}For any two elements $A_{1}$ and $A_{2}$ of algebra $A$,
the following equality holds:
\begin{equation}
\sum\limits_{k=0}^{n}\dbinom{n}{k}(A_{1}-I)^{k}(A_{2}+I)^{n-k}=\sum\limits%
_{k=0}^{n}\dbinom{n}{k}A_{1}^{k}A_{2}^{n-k}.  \label{(5)}
\end{equation}%
\end{lemma}
If $A_{1}$ and $A_{2}$ commute, then both sums in (\ref{(5)}) are equal to
$(A_{1}+A_{2})^{n}$ and thus to each other. The Lemma states that  (\ref{(5)}) still holds in the non-commutative case, when the binomial formula does not apply.
\item Another sometimes useful observation is

\begin{lemma}\label{L3} Let $\left[ D,U\right] =\lambda U$ and $U^{-1}$ exist. Then $$
B(n,\lambda ,U,D)=(DU^{-1})^{n}U^{n}.$$
\end{lemma}

\textbf{Proof} Suppose $U^{-1}$ exists. If $\left[ D,U\right] =\lambda U$,
then Lemma \ref{L1} implies that $(D-U+j\lambda I)U^{-j}=U^{-j}(D-U)$.

So, when $k=2$
$$\begin{array}{l}
(D-U)\left( D-U+\lambda I\right) U^{n-2}\\
=(D-U)U^{-1}(D-U)U^{-1}U^{n}\\
=\left( (D-U)U^{-1}\right)
^{2}U^{n}=(DU^{-1}-I)^{2}U^{n},
\end{array}
$$
when $k=3$
$$\begin{array}{l}(D-U)\left( D-U+\lambda I\right) \left( D-U+2\lambda
I\right) U^{n-3}\\
=\left( (D-U)U^{-1}\right) ^{3}U^{n}\\=(DU^{-1}-I)^{3}U^{n}
\end{array}
$$ and so on. Thus, by induction
\begin{equation}
B(n,\lambda ,U,D)=\left( \sum\limits\limits_{k=0}^{n}\dbinom{n}{k}%
(DU^{-1}-I)^{k}I^{n-k}\right) U^{n}.  \label{(4)}
\end{equation}
Using the standard binomial formula for the commuting operators $(DU^{-1}-I)$
and $I$ we get
$B(n,\lambda ,U,D)=((DU^{-1}-I)+I)^{n}U^{n}$, or $B(n,\lambda
,U,D)=(DU^{-1})^{n}U^{n}$. $\Box$

\item Note that if $U^{-1}$ exists and $[D,U]=-\lambda U,$ then the technique used
in the proof of Lemma \ref{L3} allows us to rewrite $B(n,\lambda ,U,D)$ as

$$B(n,\lambda ,U,D)=\left( \sum\limits\limits_{k=0}^{n}\dbinom{n}{k}%
(DU-U^{2})^{k}U^{2(n-k)}\right) U^{-n},$$
however the standard binomial
formula is not applicable to this sum, since the operators $DU-U^{2}$ and $%
U^{2}$ do not commute.
\end{itemize}

\bibliographystyle{amsplain}

\end{document}